\documentclass[12pt,A4paper]{article}
\usepackage[english]{babel}
\usepackage{float} 
\usepackage{array}
\usepackage[utf8]{inputenc}
\usepackage{amsthm,amsmath,amssymb,amsfonts}
\usepackage{enumerate}
\usepackage{tikz}
\usepackage{tikz-cd}
\usepackage{anysize}
\usepackage{indentfirst}
\usepackage{xypic}
\usepackage{mathabx}
\usepackage{xy}
\usepackage{mathdots}
\usetikzlibrary{arrows}

\theoremstyle{theorem}
\newtheorem{theorem}{Theorem}

\newtheorem{theo}{Theorem}[section]
\newtheorem{prop}[theo]{Proposition}

\newtheorem{cor}[theo]{Corollary}
\newtheorem{corollary}[theorem]{Corollary}
\theoremstyle{definition}
\newtheorem{df}[theo]{Definition}
\newtheorem{exa}[theo]{Example}
\newtheorem{obs}[theo]{Remark}
\newtheorem{para}[theo]{}

\newcommand{\G}{\Gamma}
\newcommand{\D}{\Delta}
\newcommand{\ra}{\rightarrow}

\newcommand{\li}{\varprojlim}

\newcommand{\N}{\mathbb{N}}
\newcommand{\Z}{\mathbb{Z}}

\newcommand{\GA}{\mathcal{G}}
\newcommand{\HA}{\mathcal{H}}

\newcommand{\KA}{\mathcal{K}}
\newcommand{\FGA}{\Pi^{abs}}

\newcommand{\F}{\mathcal{F}}

\newcommand{\bigcupdot}{\hspace{9pt}\cdot \hspace{-9pt}\bigcup}

\title{Generalized Stallings' decomposition theorem for pro-$p$ groups}
\author{Mattheus P. S. Aguiar and Pavel A. Zalesskii}

\usetikzlibrary{shapes,decorations,arrows,calc,arrows.meta,fit,positioning}
\tikzset{
	-Latex,auto,node distance =1 cm and 1 cm,semithick,
	state/.style ={ellipse, draw, minimum width = 0.7 cm},
	point/.style = {circle, draw, inner sep=0.04cm,fill,node contents={}},
	bidirected/.style={Latex-Latex,dashed},
	el/.style = {inner sep=2pt, align=left, sloped}
}

\begin{document}
	\pretolerance10000
	\maketitle
	\begin{abstract} 
		The celebrated Stallings' decomposition theorem states that the splitting of a finite index subgroup $H$ of a finitely generated group $G$  as  an amalgamated free product or  an HNN-extension over a finite group implies the same  for $G$. We generalize the pro-$p$ version of it proved by Weigel and the second author in \cite{Zalesski} to splittings over infinite pro-$p$ groups. This generalization does not have any abstract analogs. We also prove that generalized  accessibility of finitely generated pro-$p$ groups is closed for commensurability.   
	\end{abstract}



	\section{Introduction}
	
	In 1965, J-P.~Serre showed that a torsion free virtually free pro-$p$ group must be free
	(cf. \cite{S 65}). This motivated him to ask the question whether the same statement holds
	also in the discrete context. His question was answered positively some years later.
	In several papers (cf. \cite{stallbull}, \cite{stallends}, \cite{swancoh1}),
	J.R.~Stallings and R.G.~ Swan showed that free groups are precisely the groups of cohomological dimension $1$, and at the same time
	J-P.~Serre himself showed that, in a torsion free group $G$, the cohomological dimension of a
	subgroup of finite index coincides with the cohomological dimension of $G$ (cf. \cite{sercohdis}).
	
	One of the major tools for obtaining this type of result - the theory of ends - provided deep results also in the presence of torsion. The	first result to be mentioned is `Stallings' decomposition theorem' (cf. \cite{stalldeco2}). It generalizes the previously mentioned result as follows.
	
	\begin{theorem}[J.R.~Stallings]
		\label{thmstall}
		Let $G$ be a finitely generated group containing a subgroup $H$ of finite index
		which splits as a non-trivial free amalgamated product or HNN-extension over a finite group. Then $G$ also splits either as a free product with amalgamation or as an HNN-extension over a finite group.
	\end{theorem}	
	
	The pro-$p$ version of Theorem \ref{thmstall} was proved by Thomas Weigel and the second author in \cite{Zalesski} generalizing the result of W.N. Herfort and the second author in \cite{HZ}, where it was proved for virtually free pro-$p$ groups.   The objective of this paper is to show that, in the category of pro-$p$ groups, splitting theorems hold beyond splitting over finite groups. More precisely, the result  holds for a splitting over a general pro-$p$ group $K$ provided that the factors (resp. the base group) are {\it indecomposable} over any conjugate of any subgroup of $K$, i.e. do not split as a free amalgamated pro-$p$ product or pro-$p$ HNN-extension. Note that, in the pro-$p$ case, an amalgamated free pro-$p$ product or HNN-extension might be not {\it proper} (see Subsections 2.4 and 2.5), i.e. the free factors (resp. the base group) do not embed in general in the free amalgamated product (in the HNN-extension); the verb {\it split} will mean in this paper that these embeddings hold.
	
	\begin{theorem}\label{amalgam} Let $H=H_1\amalg_K H_2$ be a free amalgamated pro-$p$ product of finitely generated pro-$p$ groups $H_1,H_2$ that are indecomposable over  any conjugate of any subgroup of $K$. Let $G$ be a pro-$p$ group having $H$ as an open normal subgroup. Then $G$ splits as a free amalgamated pro-$p$ product $G=G_1\amalg_L G_2$ such that $G_i\cap H$ are contained in some conjugate of $H_i$, $i=1,2$ and $L\cap H$ is contained in some conjugate of $K$. 		
	\end{theorem}
	
	Of course, if $H_1,H_2$ do not split as a free amalgamated pro-$p$ product or HNN-extension at all (such groups called FA-groups by Serre \cite{Serre} and so we are going to use this term in the pro-$p$ context) then Theorem \ref{amalgam} holds independently of $K$.  
	
	The class of FA pro-$p$ groups is quite large and includes many important examples. All Fab pro-$p$ groups, i.e., pro-$p$ groups whose open subgroups have finite abelianization are FA pro-$p$ groups. Note that Fab pro-$p$ groups include all just-infinite pro-$p$ groups  and  play very important role in the class field theory (in particular have importance to the Fontaine-Mazur Conjecture, cf. \cite{boston}), p-adic representation theory \cite{jaikin} and include for example all open pro-$p$ subgroups of $SL_n(\Z_p)$.  The pro-$p$ completion of Grigorchuk, Gupta-Sidki groups and other branch groups are FA pro-$p$ groups as well as the Nottingham pro-$p$ group. Splittings  as  amalgamated free products of Fab analytic pro-$p$ groups occur naturally in the study of generalized RAAG pro-$p$ groups \cite[Subsection 5.5]{QSV} where it is also proved that an amalgamated free pro-$p$ product of uniformly powerful pro-$p$ groups is always proper. Thus Theorem  \ref{amalgam} applies to these splittings of generalized RAAG pro-$p$ groups.
	
	In fact, if $H_1,H_2$ are FA, we even do not need hypothesis of normality on $H$ for odd $p$.
	
	\begin{corollary} \label{amalgam absolute} Let $p>2$ and $H=H_1\amalg_K H_2$ be a free amalgamated pro-$p$ product of finitely generated FA pro-$p$ groups $H_1,H_2$. Let $G$ be a pro-$p$ group having $H$ as an open subgroup. Then $G$ splits as a free amalgamated pro-$p$ product $G=G_1\amalg_L G_2$ such that $G_i\cap H$ are contained in some conjugate of $H_i$, $i=1,2$ and $L\cap H$ is contained in some conjugate of $K$.

	\end{corollary}
	
	For an HNN-extension the corresponding statement admits two types of splittings.
	
	\begin{theorem}\label{HNN}  Let $H=HNN(H_1, K, t)$ be a  pro-$p$ HNN-extension  of a finitely generated pro-$p$ group $H_1$ that is indecomposable over  any conjugate of any subgroup of $K$. Let $G$ be a pro-$p$ group having $H$ as an open normal subgroup. Then $G$ splits as a free amalgamated pro-$p$ product $G=G_1\amalg_L G_2$  or HNN-extension $G=(G_1, L, t)$ such that  $G_i\cap H$, $i=1,2$ are contained in some conjugate of $H_1$, and $L\cap H$ is contained in some conjugate of $K$. 
		
	\end{theorem}
	
	If $H_1$ is FA then for $p>2$ we can drop the hypothesis of normality on $H$.
	
	\begin{corollary}\label{HNN not normal} Let $p>2$ and   $H=HNN(H_1, K, t)$ be a  pro-$p$ HNN-extension  of a finitely generated  FA pro-$p$ group $H_1$. Let $G$ be a pro-$p$ group having $H$ as an open subgroup. Then $G$ splits as a pro-$p$ HNN-extension $G=(G_1, L, t)$ such that  $G_1\cap H$ is contained in some conjugate of $H_1$, and $L\cap H$ is contained in some conjugate of $K$. 
		
	\end{corollary}
	
	Of course, in general, the factors of an amalgamated free pro-$p$ product $H=H_1\amalg_K H_2$ or the base group of a pro-$p$  HNN-extension $H=HNN(H_1,K,t)$ can split further, so to extend our results to a more general context we need to have some pro-$p$ version of JSJ-decomposition, i.e.  $H$ should be the fundamental pro-$p$ group of a graph of pro-$p$ groups whose  vertex groups do not split further over edge groups. Thus we need to exploit a pro-$p$ version of the Bass-Serre theory of groups acting on trees.

	\begin{theorem}\label{a virtual action} Let $G$ be a finitely generated  pro-$p$ group having an open normal  subgroup $H$ acting  on a pro-$p$ tree $T$. Suppose   $\{H_v \mid v \in V(T)\}$ is $G$-invariant. 
		Then $G$ is the fundamental group of a profinite graph of pro-$p$ groups such that each vertex  group intersected with $H$ stabilizes a  vertex of $T$. In particular $G$ splits as a non-trivial free amalgamated pro-$p$ product or a pro-$p$ HNN-extension. 
	\end{theorem} 
	
	If the stabilizers $H_v$ are FA, then the $G$-invariancy $\{H_v \mid v \in V(T)\}$ is automatic; moreover, if the $H_v$ are Fab, then  we can drop the normality assumption on $H$.
	
	\begin{corollary}\label{no normal} Let $G$ be a finitely generated  pro-$p$ group having an open  subgroup $H$ acting  on a pro-$p$ tree $T$ such that each stabilizer $H_v$ is Fab.  Then $G$ is the fundamental group of a profinite graph of pro-$p$ groups such that each vertex  group intersected with $H$ stabilizes a  vertex of $T$. In particular $G$ splits as a non-trivial free amalgamated pro-$p$ product or a pro-$p$ HNN-extension. 
	\end{corollary} 
	
	Note also that Theorem \ref{a virtual action}  does not require necessarily existence of JSJ-decomposition or even accessibility (see Section 2 for definition); in other words, we do not require that $H\backslash T$ is finite. Indeed, G. Wilkes \cite{Wilkes} constructed an example of a finitely generated  inaccessible pro-$p$ group $G$ (that acts on a pro-$p$ tree with infinite $G\backslash T$), but our theorem holds for his example as well (see Section \ref{wilkes example}). 	
	
	However, if we assume accessibility, we can tell more. 
	\begin{theorem}\label{a virtual graph of groups} Let $G$ be a finitely generated  pro-$p$ group having an open normal subgroup $H$  that splits as the fundamental pro-$p$ group of a finite graph of finitely generated pro-$p$ groups $(\HA, \Delta)$. Suppose  conjugacy classes of vertex groups are $G$-invariant. Then $G$ is the fundamental group of a reduced  finite graph of pro-$p$ $(\GA, \Gamma)$ groups such that its vertex and edge groups intersected with $H$ are subgroups of vertex and edge groups of $H$ respectively. Moreover, 	$|E(\Gamma)|\leq |E(\Delta)|$. 
	\end{theorem} 	
	
	Once more, if the vertex groups $\HA(v)$ are Fab, then we can omit $G$-invariancy and normality  hypotheses.
	
	\begin{corollary}\label{not normal} Let $G$ be a finitely generated  pro-$p$ group having an open  subgroup $H$  that splits as a finite graph of finitely generated pro-$p$ groups $(\HA, \Delta)$. Suppose  the  vertex groups of $(\HA, \Delta)$ are Fab. Then $G$ is the fundamental group of a reduced  finite graph of pro-$p$ groups $(\GA, \Gamma)$ such that its vertex and edge groups intersected with $H$ are subgroups of vertex and edge groups of $H$ respectively. Moreover, 	$|E(\Gamma)|\leq |E(\Delta)|$.
		
	\end{corollary}	
	
	\medskip
	Theorem \ref{a virtual graph of groups}  is a generalization of the pro-$p$ version of Stallings' decomposition theorem proved in \cite{Zalesski},  namely if in Theorem \ref{a virtual graph of groups} we suppose that $H$ is a non-trivial free pro-$p$ product, we obtain as a particular case 	the pro-$p$ version of  \cite[Theorem 1.1]{Zalesski}. 
	Theorem 1.4 of \cite{Wilton2} gives an example of a situation when Theorem \ref{a virtual graph of groups} is applicable, namely if all vertex groups are Poincaré duality of dimension $n$ ($PD^n$ pro-$p$ groups) and the edge groups have cohomological dimension $\leq n-1$. Moreover, many 3-manifold groups  admit a $p$-efficient JSJ-decomposition by \cite[Theorem A]{wilkes17} and if the vertex groups of them are arithmetic,   then  the pro-$p$ version of Theorem \cite[Proposition 6.23]{wilkes19} combined with \cite[Theorem 5.13]{GJPZ}  give the pro-$p$ JSJ-decomposition of their pro-$p$ completion that meets the hypothesis of Theorem \ref{a virtual graph of groups}.  
	
	The proofs of Theorem \ref{a virtual graph of groups} and Corollaries \ref{amalgam absolute}, \ref{HNN not normal}, and \ref{not normal} are more  subtle  and  require the following  theorem that is of independent interest. Note that for an open subgroup $H$ of the fundamental pro-$p$ group  $G=\Pi_1(\GA, \Gamma,v)$ of a finite graph of pro-$p$ groups, the pro-$p$ version of the Bass-Serre theorem for subgroups works, i.e. $H=\Pi_1(\HA,H\backslash S(G))$ in the standard manner (see Proposition \ref{open}).
	
	\begin{theorem}[Limitation Theorem]\label{open subgroup} Let $G=\Pi_1(\GA, \Gamma,v)$ be the fundamental pro-$p$ group of a finite  reduced graph of  pro-$p$ groups. Let $H$ be an open normal subgroup of $G$ and $H=\Pi_1(\HA,\Delta,v')$ be a decomposition as the fundamental pro-$p$ group of a reduced graph of  pro-$p$ groups $(\HA,\Delta,v')$ obtained by a reduction process from $(\HA,H\backslash S(G))$. Then $|E(\Delta)|\geq |E(\Gamma)|$.  Moreover, for $p>2$ the inequality is strict unless $\Gamma=\Delta$. 
	\end{theorem}
	
	Recall that two pro-$p$ groups $G_1, G_2$ are {\it commensurable} if there exist $H_1$ open in $G_1$ and $H_2$ open in $G_2$ such that $H_1\cong H_2$.	
	Theorem \ref{open subgroup} allows us to prove that the accessibility of a pro-$p$ group with respect to a family $\F$ of pro-$p$ groups is preserved by commensurability. For accessible abstract groups such a result can be deduced from the Stallings splitting theorem; we are not aware of such a result for accessible groups with respect to a family of infinite groups in the abstract situation. 
	
	\begin{theorem}\label{commensurability}	Let $\F$ be a family of pro-$p$ groups closed for commensurability. Let $G$ be a finitely generated pro-$p$ group and $H$ an open subgroup of $G$. Then $G$ is $\F$-accessible if and only if $H$ is $\F$-accessible. 
		
	\end{theorem}

	Note that the hypothesis of non-splitting in Theorems \ref{amalgam}, \ref{HNN}  and \ref{a virtual action} are essential. The pro-$5$ completion of the triangle group  $G=\langle x,y\mid x^5, y^5, (xy)^5\rangle$ for example contains the pro-$5$ completion $\widehat S$ of a surface group $S$ as a subgroup of index $5$. The group $\widehat S$ is a free pro-5 product of free pro-$5$ groups with cyclic amalgamation, but $G$ does not split as a non-trivial  amalgamated free pro-$5$ product. Indeed, if it does, i.e. if $G=G_1\amalg_H G_2$ then all torsion elements $x, y$ and $xy$  have to belong to some free factor up to conjugation, but then they belong to the normal closure of the same free factor, say $G_1^G$; it means that $G_1^G=G$  which is impossible, since $G/G_1^G\cong G_2/H^{G_2}\neq 1$. 
	
	\medskip	
	This paper is organized as follows. In Section $2$ we recall the elements of the pro-$p$ version of the Bass-Serre theory of groups acting on trees that will be used through the text.  Section $3$ contains the proof of the Limitation Theorem. 
	Section $4$ starts with the proof of Theorem \ref{a virtual action}. Then with Limitation Theorem in hand, we prove Theorem \ref{a virtual graph of groups}. Theorems \ref{amalgam} and \ref{HNN} then follow immediately, but their corollaries require some work. Section $5$ deals with finitely generated pro-$p$ accessible groups, where we prove Theorem \ref{commensurability}.  In the last section we show that our Theorem \ref{a virtual action} also works for Wilkes' example of a finitely generated inaccessible pro-$p$ group.
	
	\section{Main concepts of the pro-$p$ version of the Bass-Serre theory}
	
	In this section we recall the necessary notions of the pro-$p$ version of the Bass-Serre theory  (see \cite{RZ, Ribes, RZ 2000} for further details).
	
	\subsection{Pro-$p$ trees}
	
	\begin{df}[Profinite graph]
		A profinite graph is a profinite space  $\Gamma$ with a distinguished closed nonempty subset $V(\Gamma)$ called the vertex set, $E(\Gamma)=\Gamma-V(\Gamma)$ the edge set and two continuous maps $d_0,d_1:\Gamma \rightarrow V(\Gamma)$ whose restrictions to $V(\Gamma)$ are the identity map $id_{V(\Gamma)}$. We refer to $d_0$ and $d_1$ as the incidence maps of the profinite graph $\Gamma$.  
	\end{df}
	
	A morphism ($q$-morphism in the terminology of \cite{Ribes}) $\alpha:\G \ra \D$ of profinite graphs is a continuous map with $\alpha d_i=d_i \alpha$ for $i=0,1$. Note that this definition allows edges to be mapped to vertices. By \cite[Proposition 2.1.4]{Ribes} every profinite graph $\G$ is an inverse limit of finite quotient graphs of $\G$.
	
	A profinite graph is called connected if every finite quotient graph of it is connected. So a connected profinite graph is an inverse limit of finite connected graphs.
	
	\begin{df}[\cite{Mattheus}, Definition 3.4] If $\Gamma$ is a connected finite graph, its pro-$p$ fundamental group $\pi_1(\Gamma,v)$ can be defined as the pro-$p$ completion $(\pi_1^{abs}(\Gamma, v))_{\hat{p}}$ of the abstract (usual) fundamental group $\pi_1^{abs}(\Gamma, v)$. If $\Gamma$ is a connected profinite graph and $\Gamma=\li \Gamma_i$ its decomposition as inverse limit of finite graphs $\Gamma_i$, then $\pi_1(\Gamma, v)$ can be defined as the inverse limit $\pi_1(\Gamma, v)=\li (\pi_1(\Gamma_i, v_i))_{\hat{p}}$, where $v_i$ is the image of $v$ in $\Gamma_i$ (see \cite[Proposition 3.3.2 (b)]{Ribes}).  We say that $\Gamma$ is a pro-$p$ tree if $\pi_1(\Gamma)=1$. \end{df}

	If $v$ and $w$ are elements of a  pro-$p$ tree $T$, one denotes by $[v,w]$ the smallest  pro-$p$ subtree of $T$ containing $v$ and $w$.
	
	If $T$ is a pro-$p$ tree, then we say that a pro-$p$ group $G$ acts on $T$ if it acts continuously on $T$ and the action commutes with $d_0$ and $d_1$.  For $t \in V(T) \bigcupdot E(T)$ we denote by $G_t$ the stabilizer of $t$ in $G$. For a pro-$p$ group $G$ acting on a pro-$p$ tree $T$ let $\widetilde{G}$ denote the subgroup generated by all vertex stabilizers.
	
	\subsection{Fundamental pro-$p$ group of a profinite graph of pro-$p$ groups}

	\begin{df}[Sheaf of pro-$p$ groups]
		Let $T$ be a profinite space. A sheaf of pro-$p$ groups over $T$ is a triple $(\GA,\pi,T)$, where $\GA$ is a profinite space  and $\pi:\GA \ra T$ is a continuous surjection satisfying the following conditions:
		\begin{enumerate}[(a)]
			\item For every $t \in T$, the fiber $\GA(t)=\pi^{-1}(t)$ over $t$ is a pro-$p$ group (whose topology is induced by the topology of $\GA$ as the subspace topology);
			
			\begin{center}
				
			\end{center}
			\item If we define \[\GA^2=\{(g,h) \in \GA \times \GA \mid \pi(g)=\pi(h)\},\] then the map $\mu:\GA^2 \ra \GA$ given by $\mu_{\GA}(g,h)=gh^{-1}$ is continuous.
		\end{enumerate} 
	\end{df}
	
	\begin{df}
		A morphism $\underline{\alpha}=(\alpha,\alpha'):(\GA,\pi,T) \ra (\GA',\pi',T')$ of sheaves of pro-$p$ groups consists of a pair of continuous maps $\alpha: \GA \ra \GA'$ and $\alpha': T \ra T'$ such that the diagram
		\begin{equation*}
			\begin{tikzcd}
				\GA \arrow{r}{\alpha} \arrow{d}[swap]{\pi} & \GA' \arrow{d}{\pi'}\\
				T \arrow{r}[swap]{\alpha'} & T'  
			\end{tikzcd}
		\end{equation*}
		commutes and the restriction of $\alpha$ to $\GA(t)$ is a homomorphism from $\GA(t)$ into $\GA'(\alpha'(t))$, for each $t \in T$.
	\end{df}

	\begin{df}[Profinite graph of pro-$p$ groups]
		Let $\G$ be a connected profinite graph with incidence maps $d_0,d_1: \G \ra V(\G)$. A profinite graph of pro-$p$ groups over $\G$ is a sheaf $(\GA,\pi,\G)$ of pro-$p$ groups over $\G$ together with two morphisms of sheaves $(\partial_i,d_i):(\GA,\pi,\G) \ra (\GA_V,\pi,V(\G))$, where $(\GA_V,\pi,V(\G))$ is a restriction sheaf of $(\GA,\pi,\G)$ and the restriction of $\partial_i$ to $\GA_V$ is the identity map $id_{\GA_V}$, $i=0,1$; in addition, we assume that the restriction of $\partial_i$ to each fiber $\GA(m)$ is an injection.
	\end{df}
	
	\begin{obs}
		If $\Gamma$ is finite then the notion of sheaf is not needed, since $\GA=\bigcupdot_{m\in \Gamma} \GA(m)$ has the disjoint union topology. A finite graph of finite $p$-groups is just a usual graph of groups from the Bass-Serre theory. 
	\end{obs}
	
	\begin{df} A morphism of graphs of groups $\underline\nu=(\nu,\nu'):(\GA,\Gamma) \longrightarrow (\HA, \Delta)$ is a morphism of sheaves such that  $\nu\partial_i=\partial_i\nu$.
	\end{df}
	
	As was already mentioned in the introduction unlike the situation for abstract graphs of groups,  the vertex groups of a profinite graph of pro-$p$ groups $ (\GA,\Gamma)$  do not always embed in its fundamental pro-$p$ group $\Pi_1(\GA, \Gamma)$. This motivates the following definition:
	
	\begin{df}[Injective graph of pro-$p$ groups, cf. Section 6.4 of \cite{Ribes}]
		We say that a graph of pro-$p$ groups $(\GA,\G)$ is injective if the restriction of $\nu: \GA \ra \Pi_1(\GA,\G)$ to each fiber $\GA(m)$ $(m \in \G)$, is injective.
	\end{df}
	
	To achieve such embedding one has to replace the vertex and edge group with their images in $\Pi_1(\GA, \Gamma)$ (see \cite{Ribes} for details). Since we consider here splittings of pro-$p$ groups as a graph of pro-$p$ groups, we take a different  approach developed in \cite{Mattheus} to define $\Pi_1(\GA, \Gamma,v)$ with respect to a base point  that gives the embedding above automatically.  
	
	Let $I$ be a partially ordered set and $\{(\GA_i,\pi_i,\G_i),\nu_{ij}\}$ an inverse system of finite graphs of finite $p$-groups. Then $(\GA,\pi,\G)=\li_{i \in I}(\GA_i,\pi_i,\G_i)$ is a profinite graph of pro-$p$ groups. 
	
	\begin{prop}[\cite{Mattheus}, Proposition 2.15]\label{decomposition of graph of groups} Let $(\GA,\G)$ be a profinite graph of pro-$p$ groups. Then $(\GA,\G)$ decomposes as an inverse limit $(\GA,\G)=\li_{i \in I}(\GA_i,\G_i)$ of finite graphs of finite $p$-groups.
	\end{prop}	
	
	We need the following concepts to define the fundamental group of a graph of finite $p$-groups with a base point.
	
	\begin{df}[The group $F(\GA,\G)$ {\cite[Sect.\, I.5.1]{Serre}}]
		The path group $F(\GA,\G)$ is defined by $F(\GA,\G)=W_1/N$, where $W_1=\left(\mathop{\Asterisk}_{v \in V(\G)} \GA(v)\right) \Asterisk F(E(\G))$, where $F(E(\G))$ denotes the free group with basis $E(\G)$ and $N$ is a normal subgroup of $W_1$ generated by the set $\left\{ \partial_0(x)^{-1}e\partial_1(x)e^{-1} \mid x \in \GA(e), e \in E(\G) \right\}$. 
	\end{df}
	
	\begin{df}[Words of $F(\GA,\G)$ {\cite[Sec.\,I.5.1, Definition 9]{Serre}}]
		Let $c=v_0,e_0, \cdots, e_n,v_n$, be a path in $\G$ with length $n=l(c)$ such that $v_j \in V(\G), e_j \in E(\G)$, $j=0, \cdots, n$. A word of type $c$ in $F(\GA,\G)$ is a pair $(c,\mu)$ where $\mu=(g_0, \cdots, g_n)$ is a sequence of elements $g_j \in \GA(v_j)$. The element $|c,\mu|=g_0,e_0,g_1,e_1, \cdots,e_n,g_n$ of $F(\GA,\G)$ is said to be associated with the word $(c,\mu)$.   
	\end{df}	
	
	\begin{df}[The fundamental group of $(\GA,\G)$ {\cite[Sect. I.5.1, Definition 9(a)]{Serre}}]
		Let $v$ be a vertex of $\G$. We define $\pi_1(\GA,\G,v)$ as the set of elements of $F(\GA,\G)$ of the form $|c,\mu|$, where $c$ is a path whose extremities both equal $v$. One sees immediately that $\pi_1(\GA,\G,v)$ is a subgroup of $F(\GA,\G)$, called the fundamental group of $(\GA,\G)$ at $v$. In particular, if $\GA$ consists of trivial groups only then  $\pi_1(\GA,\G,v)$ becomes a usual fundamental group of the graph $\Gamma$ and denoted by  $\pi_1(\Gamma, v)$. It can be viewed of course as a subgroup that consists of set of elements of $F(\GA,\G)$ of the form $|c,\mu|=g_0,e_0,g_1,e_1, \cdots,e_n,g_n$, where $c$ is a path whose extremities both equal $v$ and $g_0=1=g_1=\ldots =g_n$. This way $G=\pi_1(\GA,\G,v)$ is a semidirect product $\pi_1(\GA,\G,v)=\langle \GA(v)\mid v\in V(\Gamma)\rangle^G \rtimes \pi_1(\Gamma, v)$.
	\end{df} 
	
	\begin{prop}[\cite{Mattheus}, Proposition 3.6] \label{fg}
		An inverse limit $(\GA,\G)=\li_{i \in I} (\GA_i,\G_i)$ of finite abstract graphs of finite $p$-groups induces an inverse limit $\li_{i \in I} (\pi_1(\GA_i,\G_i,v_{i}))_{\hat{p}}$ of the pro-$p$ completions of fundamental abstract groups $ \pi_1(\GA_i,\G_i,v_{i})$. 
	\end{prop}	
	
	\begin{df}[\cite{Mattheus}, Definition 3.7] \label{fundamental group at point} Let $(\GA,\Gamma)$ be a profinite graph of pro-$p$ groups and $(\GA,\G)=\li_{i \in I} (\GA_i,\G_i)$ be the decomposition as the inverse limit of finite  graphs of finite $p$-groups (see Proposition \ref{decomposition of graph of groups}). Let $v$ be a vertex of $\Gamma$ and $v_i$ its image in $\Gamma_i$. The group $\li_{i \in I} (\FGA_1(\GA_i,\G_i,v_{i}))_{\hat{p}}$ from Proposition \ref{fg}  will be called the pro-$p$ fundamental group of the graph of pro-$p$ groups $(\GA, \Gamma)$ at point $v$ and denoted by $\Pi_1(\GA,\G,v)$.
	\end{df}
	
	By \cite[Theorem 3.9]{Mattheus} our definition is equivalent to one in \cite{Ribes} assuming the vertex groups of $(\GA,\Gamma)$ embed in $\Pi_1(\GA, \Gamma)$, i.e., $(\GA,\Gamma)$ is injective. We denote $\nu(\GA(m))$ by $\Pi(m)$. Also note that, if $(\GA,\G)$ is a finite graph of finite $p$-groups, then $\Pi_1(\GA,\G,v)= (\FGA_1(\GA,\G,v))_{\hat{p}}$.
	
	\subsection{Reduced graph of pro-$p$ groups}
	
	\begin{df}[Reduced graph of groups] \label{rgg}
		A profinite graph of pro-$p$ groups $(\GA,\G)$ is said to be reduced if for every edge $e$, which is not a loop, neither $\partial_1:\GA(e) \ra \GA(d_1(e))$ nor $\partial_0:\GA(e) \ra \GA(d_0(e))$ is an isomorphism; we say that an edge $e$ is fictitious if it is not a loop and one of the edge maps $\partial_i$ is an isomorphism.
	\end{df}
	
	Any finite graph of groups can be transformed into a reduced finite graph of groups by collapsing fictitious edges using the following procedure. If $e$ is a fictitious edge, we can remove $\{e\}$ from the edge set of $\G$, and identify $d_0(e)$ and $d_1(e)$ to a new vertex $y$. Let $\G'$ be the finite graph given by $V(\G')=y \cup V(\G) \backslash \{d_0(e),d_1(e)\}$ and $E(\G')=E(\G) \backslash \{e\}$, and let $(\GA',\G')$ denote the finite graph of groups based on $\G'$ given by $\GA'(y)=\GA(d_1(e))$ if $\partial_0(e)$ is an isomorphism, and $\GA'(y)=\GA(d_0(e))$ if $\partial_0(e)$ is not an isomorphism. This procedure can be continued until there are no fictitious edges. The resulting finite graph of groups $(\overline{\GA},\overline{\G})$ is reduced.

	\begin{obs} \label{2.2} 
		The reduction procedure described above does not change the fundamental group (as a group given by presentation), i.e. choosing a maximal subtree to contain the collapsing edge, the morphism $(\GA, \G) \ra (\GA', \G')$ induces the identity map on the fundamental group with presentation given by eliminating redundant relations associated with fictitious edges that are just collapsed by reduction.
		
	\end{obs}		
	
	\begin{obs}	The reduction procedure can not be applied, however, if $\Gamma$ is infinite, since removal of an edge results in a non-compact object. To obtain a reduced graph of pro-$p$ groups in this case one has to reconstruct the profinite graph of pro-$p$ groups following the procedure performed in the proof of Theorem \ref{a virtual action}.
	\end{obs}
	
	The reduction procedure allows us to refine the main result of \cite{HZ} as follows:
	
	\begin{theo} \label{HZ}
		Let $G$ be a finitely generated pro-$p$ group with a free open subgroup $F$. Then $G$ is the pro-$p$ fundamental group of a reduced finite graph of finite $p$-groups $(\GA, \Gamma)$  with orders of vertex groups  bounded by $[G:F]$. Moreover, if  $G=\Pi_1(\GA', \Gamma')$ is another splitting as a reduced finite graph of finite $p$-groups then $|\Gamma|=|\Gamma'|, |V(\Gamma)|=|V(\Gamma')|, |E(\Gamma)|=|E(\Gamma')|$.
	\end{theo}	
	
	\begin{proof} By \cite[Theorem 1.1]{HZ} $G$ is the pro-$p$ fundamental group of a  finite graph of finite $p$-groups $(\GA, \Gamma)$  with orders of vertex groups  bounded by $[G:F]$ and applying the reduction procedure we get the first statement. By \cite[Theorem 7.1.2]{Ribes} maximal finite subgroups of $G$ are  exactly the vertex groups of $(\GA, \Gamma)$ and $(\GA', \Gamma')$ up to conjugation, so $|V(\Gamma)|=|V(\Gamma')|$. Now by \cite[Proposition 3.4]{CZ} $G/\widetilde G=\pi_1(\Gamma)=\pi_1(\Gamma')$ is a free pro-$p$ groups of rank $|E(\Gamma)-|V(\Gamma)|+1=|E(\Gamma')-|V(\Gamma')|+1$ implying  $|E(\Gamma)|=|E(\Gamma')|$ and $|\Gamma|=|\Gamma'|$. The proof is complete.
		
	\end{proof}
	
	Two essential particular cases of the fundamental group of a finite graph of pro-$p$ groups are amalgamated free pro-$p$ products and pro-$p$ HNN-extensions.
	
	\subsection{Free pro-$p$ products with amalgamation}

	\begin{df}[\cite{RZ}, Section 9.2]
		Let $G_1$ and $G_2$ be pro-$p$ groups and let $f_i: H \ra G_i$ $(i=1,2)$ be continuous monomorphisms of pro-$p$ groups. An amalgamated free pro-$p$ product of $G_1$ and $G_2$ with amalgamated subgroup $H$ is defined to be a pushout of $f_i$ $(i=1,2)$
		\begin{equation*}
			\begin{tikzcd}
				H \arrow{r}{f_1} \arrow{d}[swap]{f_2} & G_1 \arrow{d}{\varphi_1}\\
				G_2 \arrow{r}[swap]{\varphi_2} & G  
			\end{tikzcd}
		\end{equation*}
		in the category of pro-$p$ groups, i.e., a pro-$p$ group G together with continuous homomorphisms $\varphi_i:G_i \ra G$ $(i=1,2)$ satisfying the following universal property: for any pair of continuous homomorphisms $\psi_i:G_i \ra K$ $(i=1,2)$ into a pro-$p$ group $K$ with $\psi_1f_1=\psi_2f_2$, there exists a unique continuous homomorphism $\psi: G \ra K$ such that the following diagram is commutative:
		
		\begin{equation*}
			\begin{tikzcd}
				H \arrow{r}{f_1} \arrow{d}[swap]{f_2} & G_1 \arrow{d}[swap]{\varphi_1} \arrow[bend left]{ddr}{\psi_1}\\
				G_2 \arrow{r}{\varphi_2} \arrow[bend right]{drr}[swap]{\psi_2} & G \arrow[dotted]{dr}[swap]{\psi} \\ &&K  
			\end{tikzcd}
		\end{equation*}
		
		This amalgamated free pro-$p$ product, also referred to as free pro-$p$ product with amalgamation, is denoted by $G=G_1 \amalg_H G_2$.
	\end{df}	
	
	Following the abstract notion, we can consider $H$ as a common subgroup of $G_1$ and $G_2$ and think of $f_1$ and $f_2$ as inclusions. However, unlike the abstract case where the canonical homomorphisms \[\varphi_i^{abs}:G_i \ra G_1 \star_{H} G_2 \] $(i=1,2)$ are always monomorphisms (cf. Theorem I.1 in \cite{Serre}), the corresponding maps in the category of pro-$p$ groups \[\varphi_i:G_i \ra G_1 \amalg_{H} G_2 \] $(i=1,2)$ are not always injective. This motivates the next definition:
	
	\begin{df}
		An amalgamated free pro-$p$ product $G=G_1 \amalg_H G_2$ will be called proper if the canonical homomorphisms $\varphi_i$ $(i=1,2)$ are monomorphisms. In that case we shall identify $G_1$, $G_2$ and $H$ with their images in $G$, when no possible confusion arises.
	\end{df}
	
	Throughout the paper all free pro-$p$ products with amalgamation will be proper. 
	
	The next example shows that an  amalgamated free pro-$p$ product appears as a particular case of the fundamental pro-$p$ group of a profinite graph of pro-$p$ groups.
	
	\begin{exa}[\cite{Ribes}, Example 6.2.3(d)] \label{amalgam1}
		
		Let $G_1$, $G_2$ and $H$ be pro-$p$ groups and consider the following graph of groups:  
		
		\begin{center}
			\begin{tikzpicture}
				\node (0) at (0,0) [label=above:{\scriptsize $G_1$},point];
				\node (1) at (2,0) [label=above:{\scriptsize $G_2$},point];
				
				\draw[->] (0) edge node[below] {\scriptsize{$H$}} (1);
			\end{tikzpicture}
		\end{center}
		
		Then its fundamental pro-$p$ group will be $G=G_1 \amalg_{H} G_2$. 
	\end{exa}

	\subsection{Pro-$p$ HNN-extensions} \label{5}

	\begin{df}[\cite{RZ}, Section 9.4]
		Let $H$ be a pro-$p$ group and let $f:A \ra B$ be a continuous isomorphism between closed subgroups $A$, $B$ and $H$. A pro-$p$ $HNN$-extension of $H$ with associated groups $A$, $B$ consists of a pro-$p$ group $G=HNN(H,A,f)$, an element $t \in G$ called the stable letter, and a continuous homomorphism $\varphi: H \ra G$ with $t(\varphi(a))t^{-1}=\varphi f(a)$ and satisfying the following universal property: for any pro-$p$ group $K$, any $k \in K$ and any continuous homomorphism $\psi: H \ra K$ satisfying $k(\psi(a))k^{-1}=\psi f(a)$ for all $a \in A$, there is a continuous homomorphism $\omega: G \ra K$ with $\omega(t)=k$ such that the diagram
		\begin{equation*}
			\begin{tikzcd}
				G  \arrow[dotted]{dr}{\omega} & \\
				H \arrow{u}{\varphi} \arrow{r}[swap]{\psi} & K  
			\end{tikzcd}
		\end{equation*}
		is commutative. 
	\end{df}
	
	By construction, $G=HNN(H,A,f)$ arises as the pro-$p$ completion of the abstract $HNN$-extension $G^{abs}=HNN(H,A,f)$ (cf. \cite[Proposition 9.4.1]{RZ}). 
	
	
	In contrast with the abstract situation, the canonical homomorphism $\varphi:H \ra G=HNN(H,A,f)$ is not always a monomorphism. When $\varphi$ is a monomorphism, we shall call $G=HNN(H,A,f)$ a proper pro-$p$ $HNN$-extension. Throughout the paper all pro-$p$ $HNN$-extensions will be proper.  
	
	The next example shows that a pro-$p$ $HNN$-extension  appear as a very particular case of the pro-$p$ fundamental group of a profinite graph of pro-$p$ groups.
	
	\begin{exa}[\cite{Ribes}, Example 6.2.3(e)] \label{HNNexa}
		
		Let $\GA(v)$ and $\GA(e)$ be pro-$p$ groups and consider the following graph of groups:  
		
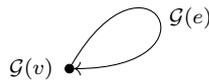
\begin{center}
			\begin{tikzpicture}[every edge/.append style={nodes={font=\scriptsize}}]
				\node (0) at (0,0) [label=left:{\scriptsize $\GA(v)$}, point];
				
				\draw[->] (0) edge[in=0,out=60,looseness=80] node[right] {$\GA(e)$} (0);
			\end{tikzpicture}
		\end{center}
		
		Then its fundamental pro-$p$ group will be  the pro-$p$ $HNN$-extension $G=\textrm{HNN}(\GA(v),\partial_0(\GA(e)),t,f)$ of $\GA(v)$, where $f:\partial_0(\GA(e)) \ra \partial_1(\GA(e))$ is the isomorphism defined by $\partial_0(x) \mapsto \partial_1(x)$, for all $x \in \GA(e)$ and $t$ is the stable letter related to $e \in E(\G)$.
	\end{exa}	
	
	\subsection{Standard pro-$p$ tree}
	\label{standard}
	Next we shall describe that standard pro-$p$ tree on which $G=\Pi_1(\GA,\G,v)$ acts. We shall assume that $\Gamma$ is finite, since we use it here only  for this case in which the notation is much simpler; see \cite[Section 6.3]{Ribes}  for the general case.  
	\begin{para}{\bf Standard (universal) pro-$p$ tree (cf. \cite{Ribes} Example 6.3.1).}
		Associated with the finite graph of pro-$p$ groups $({\cal G}, \Gamma)$ there is
		a corresponding  {\em  standard pro-$p$ tree}  (or universal covering graph)
		$S=S(G)=\bigcupdot_{m\in \Gamma}
		G/\Pi(m)$ (cf. \cite[Proposition 3.8]{ZM}).  The vertices of
		$S$ are  cosets of the form
		$g\Pi(v)$, with $v\in V(\Gamma)$
		and $g\in G$; its edges are the cosets of the form $g\Pi(e)$, with $e\in
		E(\Gamma)$; choosing a maximal subtree $D$ of $\Gamma$, the incidence maps of $S$ are given by the formulas:
		
		$$d_0 (g\Pi(e))= g\Pi(d_0(e)); \quad  d_1(g\Pi(e))=gt_e\Pi(d_1(e)) \ \ 
		(e\in E(\Gamma), t_e=1\hbox{ if }e\in D).  $$
		
		There is a natural  continuous action of
		$G$ on $S$ given by \[g(g'\Pi(m))=gg'\Pi(m),\] where $g,g' \in G$, $m \in \G$. Clearly $ G\backslash S= \Gamma$. There is a standard
		connected transversal $s:\Gamma\to S$, given by $m\mapsto \Pi(m)$.  Note
		that  $s_{|D}$  is an isomorphism of graphs and the elements $t_e$ satisfy the  equality $d_1(s(e))=t_es(d_1(e))$. Using the map $s$, we shall identify $\Pi(m)$ with the stabilizer $G_{s(m)}$ for $m\in \Gamma$:
		
		\begin{equation}\label{transversal}
			\Pi(e)=G_{s(e)}=G_{d_0(s(e))}\cap G_{d_1(s(e))}=\Pi(d_0(e))\cap t_e\Pi(d_1(e))t_e^{-1} 
		\end{equation}
		with $t_e=1$ if $e\in D$. 
		Remark also that since $\Gamma$ is finite, $E(S)$ is compact.
		
		We shall often use the following result from \cite{ZM90} stating that for open subgroups of the fundamental pro-$p$ group of finite graph of pro-$p$ groups the subgroup theorem of the Bass-Serre theory works.
		
		\begin{prop}\label{open}(\cite[Corollary 4.5 combined with 5.4]{ZM90})  Let $G=\Pi_1(\GA,\Gamma)$ be the pro-$p$ fundamental group of a finite graph of pro-$p$ groups and $H$ an open subgroup of $G$. Let $s:H\backslash S(G)\rightarrow S(G)$ be a connected transversal. Then $H=\Pi_1(\HA, H\backslash S(G))$ with $\HA(m)=H_{s(m)}$ for each $m\in H\backslash S(G)$.\end{prop}

	\end{para}

	
	
	
	
	
	\section{The proof of the Limitation Theorem for virtually free pro-$p$ groups}
	
	In this section we prove a special case of our main technical result, namely Theorem \ref{open subgroup}. We will provide the additional elements to prove Theorem \ref{open subgroup} in the next section.  
	
	\begin{theo} \label{op2} Let $G=\Pi_1(\GA, \Gamma,v)$ be the fundamental pro-$p$ group of a finite  reduced graph of  finite $p$-groups. Let $H$ be an open normal subgroup of $G$ and $H=\Pi_1(\HA,\Delta,v')$ be a decomposition as the fundamental pro-$p$ group of a reduced graph of finite $p$-groups. Then $|E(\Delta)|\geq |E(\Gamma)|$.   
	\end{theo}
	
	\begin{proof} Using induction on the index $[G:H]$ we may assume that $[G:H]=p$. Consider the action of $G$ on its standard pro-$p$ tree $S(G)$ (see Section \ref{standard}). Then  $G/H$ acts naturally on the quotient graph $H\backslash S(G)$. 
		
		Denote by $V_1$ the set of fixed vertices by this action and by $V_2$ the moved ones. By Proposition \ref{open},  $H=\Pi_1(\HA, H\backslash S(G))$ and $\HA(w)$ is a conjugate of some vertex group $\GA(v)\leq H$ for each $w\in V_2$. If $(\HA, H\backslash S(G))$ is not reduced, we can apply the procedure described after Definition \ref{rgg} to obtain the reduced graph of finite $p$-groups $(\HA, \Delta)$. Since $G$ is virtually free pro-$p$ one can use then Theorem \ref{HZ} to deduce that it suffices to prove the statement for $(\HA, \Delta)$. 
		
		Identifying $V_1$ with its bijective image in $\Gamma$ we have that for each $v\in V_1$ the vertex group  $\HA(v)=\GA(v)\cap H$ is of index $p$ in $\GA(v)$. If $V_1=\emptyset$ then all the edge and vertex groups of $(\HA, H\backslash S(G))$ are conjugates of some edge and vertex groups of $(\GA, \Gamma)$. It follows that $(\HA, H\backslash S(G))$ is reduced, since $(\GA, \Gamma)$ is by hypothesis. But  $|H\backslash E(S(G))|= p|E(\Gamma)|$ and the result follows in this case.	
		
		Assume $V_1$ is non-empty. Denote by $\Gamma(V_i)$ the spanned graph of $V_i$, $i=1,2$ and $E_{12}$ the edges that connect vertices of $V_1$ to vertices of $V_2$. If $(\HA, H\backslash S(G))$ is not reduced, then the fictitious edges can be  only the moved ones that are in  $E(\Gamma(V_1))\cup E_{12}$. Moreover,  only one such edge from its $G/H$-orbit can be collapsed. Indeed, after collapsing an edge $e \in E_{12}$ its  vertex from $V_1$ disappears (and the rest of the vertices of $Ge$ are in $\Gamma(V_2)$); on the other hand, if $e\in E(\Gamma(V_1))$ then, after collapsing it, all the other edges from its orbit become loops.  Here are the pictures for the case $p=2$, where $g \in G/H$. 
		
		\begin{figure}[H]
			\centering	
			\begin{tikzpicture}[every edge/.append style={nodes={font=\scriptsize}}]
				\node [draw, ellipse, minimum width=2.5cm, minimum height=5cm, align=center] at   (0,0)   {};
				\node [align=center] at (-2,0){$\G(V_1)$};
				
				\node [draw, ellipse, minimum width=1.5cm, minimum height=3cm, align=center] at   (3.5,2)   {};
				\node [align=center] at (5,2){$g\G(V_2)$};
				
				\node [draw, ellipse, minimum width=1.5cm, minimum height=3cm, align=center] at   
				(3.5,-2)   {};
				\node [align=center] at (5,-2){$\G(V_2)$};
				
				\node [align=center] at (1.8,-3){$E_{12}$};
				
				\node (0) at (-0.2,0.4) [label=below:{\scriptsize $v_1$}, point];
				\node (1) at (3.3,-2.4) [label=below:{\scriptsize $v_2$}, point];
				\node (2) at (3.3,1.6) [label=below:{\scriptsize $gv_2$}, point];
				\node (3) at (3.7,-1.2) [label=above:{}, point];
				\node (4) at (3.7,-1.8) [label=above:{}, point];
				\node (5) at (3.7,2.8) [label=above:{}, point];
				\node (6) at (3.7,2.2) [label=above:{}, point];
				\node (7) at (-0.1,1.3) [label=below:{}, point];
				\node (8) at (0.5,1.3) [label=below:{}, point];
				
				\draw[->] (0) edge node[right] {\scriptsize{$e_{12}$}} (1);
				\draw[->] (0) edge node[right] {\scriptsize{}} (2);
				\draw[->] (3) edge[bend right=30] node[right] {\scriptsize{}} (4);
				\draw[->] (4) edge[bend right=30] node[right] {\scriptsize{}} (3);
				\draw[->] (5) edge[bend right=30] node[right] {\scriptsize{}} (6);
				\draw[->] (6) edge[bend right=30] node[right] {\scriptsize{}} (5);
				\draw[->] (7) edge[bend right=30] node[below] {\scriptsize{$e_1$}} (8);
				\draw[->] (8) edge[bend right=30] node[below] {\scriptsize{}} (7);
				\draw[->] (8) edge node[below] {\scriptsize{}} (5);
				\draw[->] (8) edge node[below] {\scriptsize{}} (3);
			\end{tikzpicture}
			\caption{Graph of groups $(\HA,S(G)/H)$}
		\end{figure}

		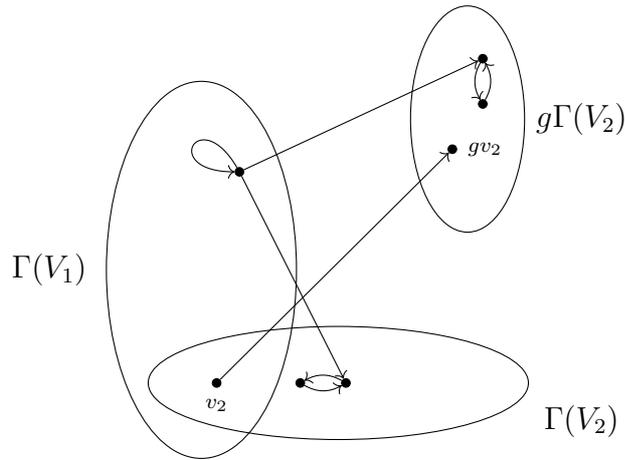
\begin{figure}[H]
			\centering	
			\begin{tikzpicture}[every edge/.append style={nodes={font=\scriptsize}}]
				\node [draw, ellipse, minimum width=2.5cm, minimum height=5cm, align=center] at   (0,0)   {};
				\node [align=center] at (-2,0){$\G(V_1)$};
				
				\node [draw, ellipse, minimum width=1.5cm, minimum height=3cm, align=center] at   (3.5,2)   {};
				\node [align=center] at (5,2){$g\G(V_2)$};
				
				\node [draw, ellipse, minimum width=5cm, minimum height=1.5cm, align=center] at   
				(1.8,-1.5)   {};
				\node [align=center] at (5,-2){$\G(V_2)$};
				
				\node (0) at (0.2,-1.5) [label=below:{\scriptsize $v_2$}, point];
				\node (2) at (3.3,1.6) [label=right:{\scriptsize $gv_2$}, point];
				\node (3) at (1.9,-1.5) [label=above:{}, point];
				\node (4) at (1.3,-1.5) [label=above:{}, point];
				\node (5) at (3.7,2.8) [label=above:{}, point];
				\node (6) at (3.7,2.2) [label=above:{}, point];
				\node (8) at (0.5,1.3) [label=below:{}, point];

				\draw[->] (8) edge[in=180,out=120,looseness=40] node[left] {} (8);
				\draw[->] (0) edge node[right] {\scriptsize{}} (2);
				\draw[->] (3) edge[bend right=30] node[right] {\scriptsize{}} (4);
				\draw[->] (4) edge[bend right=30] node[right] {\scriptsize{}} (3);
				\draw[->] (5) edge[bend right=30] node[right] {\scriptsize{}} (6);
				\draw[->] (6) edge[bend right=30] node[right] {\scriptsize{}} (5);
				\draw[->] (8) edge node[below] {\scriptsize{}} (5);
				\draw[->] (8) edge node[below] {\scriptsize{}} (3);
			\end{tikzpicture}
			\caption{Reduced graph of groups $(\HA,\D)$ assuming $e_1$ and $e_{12}$ are collapsed}
		\end{figure}
		
		Thus we can deduce that $E(\Delta)\geq |E(\Gamma(V_1))| + p|E(\Gamma(V_2))|+ (p-1)|E_{12}|\geq |E(\Gamma(V_1)| + |E(\Gamma(V_2))|+ |E_{12}|\geq  E(\Gamma)$. 
		This finishes the induction and concludes the theorem.
	\end{proof}

	\bigskip
	
	\begin{obs}\label{reverse inequality} 
		
		It follows from the first 3 lines of the proof of Theorem \ref{op2} that $|\Delta|\leq [G:H] |\Gamma|$, $|E(\Delta)|\leq [G:H] |E(\Gamma)|$ and $|V(\Delta)|\leq [G:H] |V(\Gamma)|$. 
	\end{obs}

	\begin{obs}\label{p>2} The proof also shows that for $p>2$ one has $|E(\Gamma)|< |E(\Delta)|$ unless $\Gamma=\Delta$. This means that if $[G:H]> |E(\Delta)|$ then there exists intermediate subgroup $H\leq K < G$ such that $K=\Pi_1(\KA, \Gamma,v)$. For $p=2$ the equality $|E(\Delta)|=|E(\Gamma)|$ can happen either in the case $\Delta=\Gamma$ or if for every edge $e\in E_{12}$ one has $[\GA(d_0(e)):\GA(e)]=2$.
		
	\end{obs}

	In fact Remark \ref{p>2} combined with \cite[Lemma 3.2]{SnZ} gives a short proof of \cite[Theorem B]{Kochloukova2}. We refer to the latter paper for the definition of a pro-$p$ limit group.
	
	\begin{cor} Let $G$ be a pro-$p$ limit group and $U$ be a proper open subgroup of $G$. Then the minimal number of generators $d(U)$ is strictly bigger than the minimal number of generators $d(G)$.
	\end{cor} 
	
	\section{Finitely generated pro-p groups virtually acting on trees}	
	
	In this section we prove the main results stated in the introduction and deduce several consequences. The proof of Theorem \ref{a virtual action}  follows the proof of \cite[Lemma 4.1]{CZ} whose original idea appears in the proof of the main result of \cite{Zalesski}.	
	
	\bigskip 
	
	{\it Proof of Theorem \ref{a virtual action}.} 
	
	Let $\mathcal{U}$ be the collection of open normal subgroups $U$ of $G$ contained in $H$.  Denote by $\widetilde{U}$ the topological closure of $U$ generated by the $U$-stabilizers of the vertices of $T$, i.e., \[\widetilde{U}=cl(\langle U \cap H_v \mid v \in V(T)\rangle).\] Then $\widetilde{U}$ is a closed normal subgroup of $G$. To see this it suffices to observe that  $H^{g}_v\leq \tilde H$ for any $g\in G$ which is exactly our hypothesis. 
	
	Note that $\widetilde{U} \backslash T$ is a pro-$p$ tree and  $H/\widetilde{U}$ acts on $\widetilde{U} \backslash T$ with $U/\tilde U$ acting freely. Therefore   $G/\widetilde{U}$ contains the open normal subgroup $U/\widetilde{U}$ which is finitely generated and free pro-$p$ (cf. \cite[Theorem 2.6]{ZM}). By Theorem \ref{HZ}, $G/\widetilde{U}$ is isomorphic to the pro-$p$ fundamental group $\Pi_1(\GA_U,\G_U,v_U)$ of a finite graph of finite $p$-groups. 
	
	Although neither the finite graph $\G_U$ nor the finite graph of finite $p$-groups $\GA_U$ are uniquely determined by $U$ (resp. $\widetilde{U}$), the index $U$ in the notation shall express that both these objects are depending on $U$. Using the procedure described after Definition \ref{rgg} we have a morphism $\eta:(\GA_U,\G_U) \ra (\overline{\GA}_U,\overline{\G}_U)$ to a reduced graph of groups. 
	
	For $V \subseteq U$ both open and normal in $G$ the decomposition $G/\widetilde{V}=\Pi_1(\overline{\GA}_V,\overline{\G}_V,\overline{v}_V)$ gives rise to a natural decomposition of $G/\widetilde{U}$ as the fundamental group $G/\widetilde{U}=\Pi_1(\GA_{V,U},\overline{\G}_V,\overline{v}_V)$ of a finite graph of finite $p$-groups $(\GA_{V,U},\overline{\G}_V)$, where the vertex and edge groups are $\GA_{V,U}(x)=\overline{\GA}_V(x)\tilde{U}/\tilde{U}$, $x \in \overline{\G}_V$. Thus we have a morphism $\nu_{V,U}:(\overline{\GA}_V,\overline{\G}_V) \ra (\GA_{V,U},\overline{\G}_V)$ of graphs of groups such that the induced homomorphism on the pro-$p$ fundamental groups coincides with the canonical projection $\varphi_{V,U}: G/\widetilde{V} \ra G/\widetilde{U}$. Choose a reduction morphism $\eta_U:(\GA_{V,U},\overline{\G}_V) \ra (\overline{\GA}_{V,U},\overline{\G}_U)$ to a reduced graph of groups $(\overline{\GA}_{V,U},\overline{\G}_U)$ (it is not unique); it induces the identity map on the fundamental group $G/\widetilde{U}$ (see Remark \ref{2.2}) and so $\eta_U\nu_{V,U}$ induces the homomorphism $\Pi_1(\overline{\GA}_V,\overline{\G}_V,\overline{v}_V) \ra \Pi_1(\overline{\GA}_{V,U},\overline{\G}_U,\overline{v}_U)$ on the pro-$p$ fundamental groups that coincides with the canonical projection $\varphi_{UV}: G/\widetilde{V} \ra G/\widetilde{U}$. 
	
	Using the aforementioned reduction, we have that $G/\widetilde{U}=\Pi_1(\overline{\GA}_U,\overline{\G}_U,\overline {v}_U)$. Then, by \cite[Corollary 3.3]{Zalesski}, the number of isomorphism classes of finite reduced graphs of finite $p$-groups $(\GA'_U,\G')$ which are based on $\G'$ satisfying $G/\widetilde{U} \simeq \Pi_1(\GA'_U,\G',v_0)$ is finite. 
	
	Suppose that $\Omega_U$ is a set containing a copy of every such isomorphism class. Since $G$ is finitely generated, we may choose $V_i$, $i \in \mathbb{N}$, be a decreasing chain of open normal subgroups of $G$ with $V_0 = U$ and $\bigcap_i V_i=\{1\}$. For $X \subseteq \Omega_{V_{i}}$ define $T(X)$ to be the set of all reduced graphs of groups in $\Omega_{V_{i-1}}$ that can be obtained from graphs of groups of $X$ by the procedure of reduction explained above (note that $T$ is not a map). Define $\Omega_1 = T(\Omega_{V_1})$, $\Omega_2 = T(T(\Omega_{V_2}))$, $\cdots$, $\Omega_i = T^i(\Omega_{V_i})$ and note that it is a non-empty subset of $\Omega_U$ for every $i \in \mathbb{N}$. Clearly $\Omega_{i+1} \subseteq \Omega_i$ and since $\Omega_U$ is finite there is an $i_1 \in \mathbb{N}$ such that $\Omega_j = \Omega_{i_1}$ for all $j > i_1$ and we denote this $\Omega_{i_1}$ by $\Sigma_U$. Then $T(\Sigma_{V_i}) = \Sigma_{V_{i-1}}$ and so we can construct an infinite sequence of graphs of groups $(\GA_{V_j},\G_j) \in  \Omega_{V_j}$ such that $(\GA_{V_{j-1}},\G_{j-1}) \in  T(\GA_{V_j},\G_j)$ for all	$j$. This means that $(\GA_{V_j},\G_{V_j})$ can be reduced to $(\GA_{V_{j-1}},\G_{j-1})$, i.e. this sequence $\{(\GA_{V_j},\G_j)\}$ is an inverse system of reduced graph of groups satisfying the required conditions. Therefore   $(\GA,\G)=\li (\GA_{V_j}, \Gamma_j)$ is a reduced finite graph of finitely generated pro-$p$ groups satisfying $G \simeq \Pi_1(\GA,\G,v)$.
	
	Moreover, denoting by $x_V$ the image of $x\in \Gamma$ in $\Gamma_V$ we have  $\GA(x)=\li \GA_{V_j}(x_{V_j})$ if $x$ is either a vertex or an edge of $\G$. Since $\GA_{V_j}(x)\cap H/\widetilde V_j$ fixes a vertex in $\widetilde V_j\backslash T$ for each $V_j$, and the set of fixed vertices of $\GA_{V_j}(x)\cap H/\widetilde V_j$ is compact, the projective limit argument implies that $\GA(x)\cap H$ fixes a vertex of $T$. 
	
	By \cite[Theorem 4.2]{CZ} a finitely generated pro-$p$ group that acts on a pro-$p$ tree splits as an amalgamated free pro-$p$ product  or pro-$p$ HNN-extension over the stabilizer of an edge. Using the fact that the fundamental pro-$p$ group of a graph of pro-$p$ groups acts on its standard pro-$p$ tree (see \cite[Chapter 6]{Ribes}) we can deduce that $G$ splits as non-trivial free amalgamated pro-$p$ product or pro-$p$ HNN-extension. This finishes the proof of the theorem.
	
	\hspace{14.9cm} \qedsymbol
	\bigskip
	
	\bigskip
	{\it Proof of Corollary \ref{no normal}.} 
	
	Let $N=H_G$ be the normal core of $H$ in $G$. Since $H_v$ is Fab, so is and $N_v$ and therefore  $N_v^g$ must fix a vertex of $T$.   Hence hypotheses of Theorem \ref{a virtual action} are satisfied for $N$ and the result follows.  
	
	\hspace{14.9cm} \qedsymbol
	\bigskip
	
	\begin{cor}\label{limitation}  $|E(\Gamma)|\leq  |E(H\backslash T)|$. Moreover, if $p>2$ and $\Gamma\neq H\backslash T$, then the inequality is strict.
	\end{cor}
	
	\begin{proof} It make sense to prove the statement assuming that $H\backslash T$ is finite.
		By Theorem \ref{op2} combined with Theorem \ref{HZ}, $|E(\G_{V_j})|\leq |E(H\backslash T)|$. Hence $|E(\Gamma)|\leq |E(H\backslash T)|$ as required.  Moreover, if $\Gamma\neq H\backslash T$ then 	$|E(\G_{V_j})|< |E(H\backslash T)|$ by Remark \ref{p>2} and hence $|E(\Gamma)|< |E(H\backslash T)|$.
	\end{proof}	
	
	We are ready to prove Theorem \ref{open subgroup}.
	It will be crucial to deal with the generalized version of Stallings' decomposition theorem and with accessibility in finitely generated pro-$p$ groups.
	
	\bigskip
	{\it Proof of Theorem \ref{open subgroup}}. In this case we consider the action of $H$ on the standard pro-$p$ tree $S(G)$ of $G$ and so $\widetilde U$ from the proof of Theorem \ref{a virtual action} is automatically normal in $G$. Moreover, all $\Gamma_j$ in the proof of Theorem \ref{a virtual action} can be assumed to be equal to $\Gamma$. Hence Corollary \ref{limitation} valid for this case. This finishes the proof.   
	
	\begin{obs} Corollary \ref{limitation} shows that if  $[G:H]> |E(H\backslash T)|$ then there exists an intermediate open subgroup $H\leq K\leq G$ such that $K$ splits as the fundamental pro-$p$ group of a graph of pro-$p$ groups with the same underlying graph. To illustrate this, let $H=H_1\amalg_L H_2$ be a (proper) free amalgamated product of $FA$ pro-$p$ groups (for example open pro-$p$ subgroups of $SL_n(\Z_p)$) and $G$ be a finite extension of $H$, i.e. $G/H=P$ with $P$ being finite $p$-group. Then $G= G_1\amalg_M G_2$.   Moreover,  $G/H=G_1/(G_1\cap H)\amalg_{M\cap H} G_2/(G_2\cap H)$ has to be fictitious, which means that, up to renumbering, $P=G_1/(H\cap G_1)$ and $G_2/(G_2\cap H)=M/(M\cap H)$ since $G/H$ is finite. 	
	\end{obs}	
	
	One of the obstacles to obtain the main structure result in the pro-$p$ version of Bass-Serre theory is  that a maximal subtree of a profinite graph $\Gamma$ does not always exist. The next corollary shows that, for the finitely generated case, this difficulty can be surpassed.
	
	\begin{cor} $\Gamma$ posseses a closed maximal pro-$p$ subtree.
	\end{cor}
	
	\begin{proof} By \cite[Section 2.3, Corollary 2]{Serre}  the inverse image of a maximal subtree under a collapse is a maximal subtree. Hence we can choose maximal subtrees $D_j$ of $\Gamma_{V_j}$  from the proof of Theorem \ref{a virtual action} such that they form the inverse subsystem. Then $D= \varprojlim D_j$ is a pro-$p$ tree with $V(D)=V(\Gamma)$. 
	\end{proof}  
	
	\bigskip
	
	{\it Proof of Theorem \ref{a virtual graph of groups}.} 
	
	Consider the action of $H$ on its standard pro-$p$ tree $S(H)$ (see Section \ref{standard}) and apply Theorem \ref{a virtual action} together with Corollary \ref{limitation} to get a splitting of $G$ as the fundamental group $\Pi_1(\GA, \Gamma, v)$ of a finite graph of pro-$p$ groups, where the vertex groups intersected with $H$ are subgroups of vertex groups of $H$. Now observe that Corollary \ref{limitation} or Theorem \ref{a virtual action} allows us to assume that all $\Gamma_{V_j}=\Gamma$ in the proof of Theorem \ref{a virtual action} and applying it we get the last statement of the theorem.
	
	We are left with the statement about edge stabilizers. Assume w.l.o.g. that $(\HA, \Delta)$ is reduced. Then  the vertex stabilizers of $\widetilde V_j\backslash S(H)$ in $H/\widetilde V_j$ are exactly  maximal finite subgroups of the  $H/\widetilde V_j$. This implies, in particular, that, for $V_{j+1}\leq_o V_j$, the maximal finite subgroups of $H/\widetilde V_{j+1}$ map onto maximal finite subgroups of $H/\widetilde V_j$. It induces the bijection of the conjugacy classes of the maximal finite subgroups of $H/\widetilde V_{j+1}$ and $H/\widetilde V_j$. Then, if $e$ is an edge of $\Gamma$, starting from some $j$, one has $\GA_{V_j}(e)=\GA'_{V_j}(v)\cap \GA_{V_j}(w)^g$, where $v,w$ are the extremity vertices of $e$ and are maximal finite subgroups of $G/\widetilde V_j$. 
	
	Let $H_{V_j}=\Pi_1(\HA_j, \Delta_j, v_j')$ be the splitting of $H_{V_j}$ as an open subgroup of $\Pi_1(\GA_{V_j}, \Gamma_j,v_j)$ (see Proposition \ref{open}). Then $H_{V_j}\cap \GA_{V_j}(e)\leq \HA_j(e')^h$ for some $e'\in E(\Delta_j)$, $h\in H_{V_j}$. It follows that $H_{V_j}\cap \GA_{V_j}(e)$ is contained in the intersection of at least two distict maximal finite subgroups of $H_{V_j}$ (some vertex stabilizers of $\widetilde H_{V_j}\backslash S(H)$).    Hence  $\GA(e)\cap H$ is contained in the intersection of at least two distict vertex stabilizers of $T$ and so fixes an edge of $ T$. This finishes the proof.
	
	\hspace{14.9cm} \qedsymbol
	\bigskip
	
	{\it Proof of Corollary \ref{not normal}.}
	
	Let $N=H_G$ be the normal core of $H$ in $G$. Since $\HA(v)$, $v\in V(\Delta)$ is Fab,  $N_v$ must fix a vertex of the standard pro-$p$ tree $S(G)$ (see Section \ref{standard}) and so its conjugacy class is $G$-invariant.   Hence  the result follows from Theorem \ref{a virtual graph of groups}.

	\hspace{14.9cm} \qedsymbol
	\bigskip
	
	{\it Proof of Theorem \ref{amalgam}.} 
	
	By Theorem \ref{a virtual graph of groups},  $G$ is the fundamental group of a graph of pro-$p$ groups with one edge only (cf. Example \ref{amalgam1}). However, $H/\tilde H$ is trivial in this case and since $\tilde H\leq \tilde G$, $G/\tilde G$ is finite. Hence $G$ can not be an HNN-extension.  
	
	\hspace{14.9cm} \qedsymbol
	\bigskip
	
	{\it Proof of Corollary} \ref{amalgam absolute}.
	We use induction on $[G:H]$.  The base of induction $[G:H]=p$ follows from Theorem \ref{amalgam} as 	$H$ is normal in $G$ in this case. Suppose $[G:H]>p$ and $H< N\leq G$ with $[N:H]=p$. Then, by Theorem \ref{amalgam}, $N=N_1\amalg_M N_2$ and $N_1\cap H, N_2\cap H$ are conjugate into $H_1$ or $H_2$. To apply the induction step, we just need to show that $N_1$ and $N_2$ are FA.
	
	Let $S(N)$ be a standard pro-$p$ tree on which $N$ acts. Since $H_1,H_2$ are FA, each of them fix a vertex in $S(N)$ and hence is conjugate into  $(N_1\cap H)$ or $(N_2\cap H)$.     Suppose w.l.o.g $H_1 \leq N_1\cap H$. Since $H_1$ is not conjugate into $H_2$  we deduced that   $N_1\cap H$ is conjugate into $H_1$ and hence is equal to $H_1$. 
	Then  $H_1$ has at most index $p$ in $N_1$ and so $N_1$ can not split, because otherwise  $H_1$ would split by Proposition \ref{open} and this splitting is non-trivial by Theorem \ref{open subgroup}, so $N_1$ is FA by \cite[Theorem 4.2]{CZ}.
	
	\smallskip	
	If $N_2\cap H$ is conjugate into $H_2$ then by the same argument one deduces that $N_2$ is FA.
	
	\smallskip	 
	We claim that $N_2\cap H$ is not conjugate into $H_1$. Suppose it is, so $H$ is contained in the normal closure $G_1^G$ in $G$ and, by Proposition \ref{open}, $G_1^G$ splits as a fundamental graph of groups $(\GA_1,\Delta)$ that we may assume to be reduced (see Remark \ref{2.2}). Moreover,  by Theorem \ref{open subgroup}, $E(\Delta)=1$ only if $p=2$. Since $p=2$  case is excluded by the hypothesis, the proof is complete.
	
	\hspace{14.9cm} \qedsymbol
	
	\bigskip	
	Theorem \ref{HNN} follows by direct application of Theorem \ref{a virtual action} and Corollary \ref{limitation} (cf. Example \ref{HNNexa}). 
	
	\bigskip
	
	{\it Proof of Corollary} \ref{HNN not normal}. By Corollary \ref{not normal} we just need to prove the last statement of the corollary. 
	
	We use induction on $[G:H]$.  The base of induction $[G:H]=p$ follows from Theorem \ref{HNN} as 	$H$ is normal in $G$ in this case. Suppose $[G:H]>p$ and $H< N\leq G$ with $[N:H]=p$. Then, by Theorem \ref{amalgam}, either $N=N_1\amalg_M N_2$ with  $N_1\cap H, N_2\cap H$  conjugate into $H_1$ or $N=HNN(N_1, M, t)$ with $N_1\cap H$ is conjugate into $H_1$. 
	But for $p>2$ the first case $N=N_1\amalg_M N_2$ does not occur by Theorem \ref{open subgroup}.
	Thus to apply the induction step, we just need to show that $N_1$ is FA.
	
	Let $S(N)$ be a standard pro-$p$ tree on which $N$ acts. Since $H_1$ is FA, it fixes a vertex in $S(N)$ and hence is conjugate into  $(N_1\cap H)$. Hence $H_1$ and     
	$N_1\cap H$ are conjugate  so, w.l.o.g, we may assume that with $N_1\cap H= H_1$. Then  $H_1$ has index at most $p$ in   $N_1$. Then $N_1$ can not split, because then $H_1$ would split non-trivially  by Proposition \ref{open} and Theorem \ref{open subgroup} contradicting the hypothesis, so $N_1$ is FA by \cite[Theorem 4.2]{CZ}.
	
	\hspace{14.9cm} \qedsymbol
	\bigskip
	
	\section{Generalized accessible pro-$p$ groups}
	
	Abstract accessibility was studied in a series of papers by M.J. Dunwoody (cf. \cite{D85} \cite{Dun93},\cite{DS99},\cite{DS00}), where he proved that every finitely presented group is accessible, but not every finitely generated group over an arbitrary family of groups. In fact, he presented an example of a finitely generated inaccessible group. Generalized accessible groups were studied by Bestvina and Feighn 
	(\cite{BF}). The pro-$p$ version of accessibility was introduced by G. Wilkes in \cite{Wilkes}, and  Z. Chatzidakis and the second author generalized this definition as follows:
	
	\begin{df}[Generalized accessible pro-p group, cf. Definition 5.1 of \cite{CZ}]
		Let $\F$ be a family of pro-$p$ groups. We say that a pro-$p$ group $H$ is $\F$-accessible if any splitting of $H$  as the fundamental group of a reduced finite graph $(\GA,\Gamma)$ of pro-$p$ groups such that the edge groups are in $\F$   has bound on $\Gamma$.
	\end{df}
	
	Now we prove Theorem \ref{commensurability}.
	
	\bigskip
	
	{\it Proof of Theorem \ref{commensurability}.} Using the obvious induction on $[G:H]$ we may assume that $[G:H]=p$ and so $H$ is normal in $G$.
	
	Suppose $H$ is $\F$-accessible and  $G$ is not. Then for any $n\in \N$ there exists a finite reduced graph of pro-$p$ groups  $(\GA,\Gamma)$ such that $G=\Pi_1(\GA,\Gamma,v)$ with edge groups in $\F$ and $|E(\Gamma)|>n$.  It follows  from the proof of Theorem \ref{a virtual action} that there exists an open normal subgroup $U$ of $G$ contained in $H$ such that  $G/\widetilde U=(\GA_U, \Gamma)$ is the fundamental group of a reduced quotient graph  of finite $p$-groups of $(\GA_U, \Gamma)$ over the same underlying graph $\Gamma$.  Then, by Theorem \ref{open subgroup}, $H/\tilde U=\Pi_1(\HA_U, \Delta_U,v_U)$ is the fundamental group of a finite reduced graph of pro-$p$ groups with  $|E(\Delta_U)|>n$. Then $|E(\Delta_V)|>n$ for each open normal $V$ contained in $U$. By the proof of Theorem \ref{a virtual action}, the set $\{(\HA_V, \Delta_V)\mid V\leq_o U\}$ contains a subset that form a surjective inverse system $\{(\HA_{V_j}, \Delta_{V_j})\}$ with $(\HA,\Delta)=\li (\HA_{V_j}, \Delta_{V_j})$ being the reduced graph of pro-$p$ groups such that and $H=\Pi_1(\HA, \Delta)$. Moreover, it is proved in Theorem \ref{a virtual graph of groups} that edge groups of $(\HA,\Delta)$ are virtually $\F$.  Therefore, $|E(\Delta)|>n$ for an arbitrary chosen $n\in\N$ contradicting $\F$-acessibility of $H$.
	
	Suppose now $G$ is $\F$-accessible with accessibility number $m$ and $H$ is not. Then for any $n\in \N$ there exists a finite reduced graph of pro-$p$ groups  $(\HA,\Delta)$ such that $H=\Pi_1(\HA,\Delta)$ with edge groups in $\F$ and $|E(\Delta)|>n$. Again  it follows from the proof of Theorem \ref{a virtual action} that there exists an open normal subgroup $U$ of $G$ contained in $H$ such that  $H_U=(\HA_U, \Delta,v)$ is the fundamental group of a reduced quotient graph  of finite $p$ groups with th same underlying graph $\Delta$.    On the other hand, the graph of groups $(\overline \GA_U,\overline\Gamma_U)$ with $\widetilde G_U=\Pi_1(\overline\GA_U, \overline\Gamma_U,\overline{v}_U)$  constructed in  the proof of Theorem \ref{a virtual action} must have at most $m$ edges  and therefore by Theorem \ref{HZ} and Remark \ref{reverse inequality},  $\Delta$   has at most $m[G:H]$ edges. This contradiction completes the proof of the theorem.  
	
	\section{Adaptation of Wilkes' example} \label{wilkes example}
	
	In this section we show that our Theorem \ref{a virtual action} also works for the inaccessible finitely generated group presented by Wilkes in \cite[Section 4.2]{Wilkes}. 
	
	\begin{exa} First define the map $\mu_n:\{0, . . . , p^{n+1}-1\}\longrightarrow \{0, . . . , p^n-1\}$ by sending an integer to its remainder modulo $p^n$.  Define $H_n=\mathbb{F}_p[\{0, . . . , p^n-1\}]$ to be the $\mathbb{F}_p$-vector space with basis $\{h_0, . . . , h_{{p^n}-1}\}$. There are inclusions $H_n\subseteq H_{n+1}$ given by inclusions of bases, and retractions $\eta_n:H_{n+1}\longrightarrow H_n$ defined by $h_k\rightarrow h_{{\mu_n}(k)}$. Note also that there is a natural action of $\Z/p^n\Z$ on $H_n$ given by cyclic permutation of the basis elements, and that these actions are compatible with the retractions $\eta_i$. The inverse limit of the $H_n$ along these retractions is the completed group ring $H_{\infty}=\mathbb{F}_p[[\Z_p]]$ with multiplication ignored. The continuous action of $\Z_p$ on the given basis of $H_{\infty} $ allows  to form a sort of a pro-$p$ wreath product $H_\omega=\mathbb{F}_p[[\Z_p]]\rtimes \Z_p= \varprojlim (H_n\rtimes \Z/p^n)$ which is a pro-$p$ group into which $H_{\infty}$ embeds.
		
		Next set $K_n=\mathbb{F}_p\times H_n=\langle k_n\rangle\times H_n$. Set $G_1=K_1\times\mathbb{F}_p$. For $n >1$, let $G_n$ be a finite $p$-group with presentation $G_n=\langle k_{n-1}, k_n,h_0, \ldots , h_{p^n-1} \mid k_i^p=h_i^p= 1, h_i\leftrightarrow h_j,k_{n-1}\leftrightarrow h_i\  {\rm for\ all\ } i\neq p^{n-1}$, $k_n= [k_{n-1}, h_{p^{n-1}}]\  {\rm central}\rangle$ where $\leftrightarrow$ denotes the relation ‘commutes with’. 
		
		The choice of generator names describes maps $H_n\longrightarrow G_n$, $K_{n-1}\longrightarrow G_n$, and $K_n\longrightarrow G_n$.  One may easily see that all three of these maps are injections. Define a retraction map $$\rho_n:G_n\longrightarrow K_{n-1}$$ by killing $k_n$ and by sending $h_k\rightarrow h_{{\mu_{n-1}}(k)}$. Note that $\rho_n$ is compatible with $\eta_n:H_n\rightarrow H_{n-1}$ that is, there is a commuting diagram
		
		\begin{equation*}
			\begin{tikzcd}
				K_{n-1} \arrow[hook, shift left]{r} & G_n \arrow[shift left]{l}{\rho_n}\\
				H_{n-1} \arrow[hook]{u} \arrow[hook, shift left]{r} & H_n \arrow[hook]{u} \arrow[shift left]{l}{\eta_n}  
			\end{tikzcd}
		\end{equation*}
		
		Define $\Pi_1(\GA_m,\Gamma_m,v_m)$ to be the pro-$p$ fundamental group of the following graph of groups:
		
		\begin{center}
			\begin{tikzpicture}[every edge/.append style={nodes={font=\scriptsize}}]
				\node (0) at (0,0) [label=above:{\scriptsize $G_1$},point];
				\node (1) at (2,0) [label=above:{\scriptsize $G_2$},point];
				\node (2) at (3,0) {$\cdots$};
				\node (3) at (4,0) [label=above:{\scriptsize $G_{m-1}$},point];
				\node (4) at (6,0) [label=above:{\scriptsize $G_m$},point];
				
				\draw[->] (0) edge node[above] {$K_1$} (1);
				\draw[->] (3) edge node[above] {$K_{m-1}$} (4);
			\end{tikzpicture}
		\end{center}
		
		Note that the retraction $\rho_n:G_n\longrightarrow K_{n-1}$ induces the retraction $P_{m+1}\longrightarrow P_m$  represented by the collapse the last right edge of the picture. 
		
		Then $P=\varprojlim_{m\in\N} \Pi_1(\GA_m,\Gamma_m,v_m)$ is the fundamental group of the following profinite graph of pro-$p$ groups
		
		\begin{center}
			\begin{tikzpicture}[every edge/.append style={nodes={font=\scriptsize}}]
				\node (0) at (0,0) [label=above:{\scriptsize $G_1$},point];
				\node (1) at (2,0) [label=above:{\scriptsize $G_2$},point];
				\node (2) at (4,0) [label=above:{\scriptsize $G_{3}$},point];
				\node (3) at (5,0) {$\cdots$};
				\node (4) at (6,0) [label=above:{\scriptsize $G_{\infty}$},point];
				
				\draw[->] (0) edge node[above] {$K_1$} (1);
				\draw[->] (1) edge node[above] {$K_{2}$} (2);
			\end{tikzpicture}
		\end{center}
		where the vertex at infinity is a one point compactification of the edge set of the graph and so does not have an incident edge to it; thus the edge set is not compact. The vertex group $G_{\infty}$ of the vertex at infinity is $G_\infty=K_{\infty}=\varprojlim_{i\in \N} K_i=H_{\infty}$. Let $J=P\amalg_{H_{\infty}} H_{\omega}$.  Then $J$ is the fundamental group of the following profinite graph of groups
		
		\begin{center}
			\begin{tikzpicture}[every edge/.append style={nodes={font=\scriptsize}}]
				\node (0) at (0,0) [label=above:{\scriptsize $G_1$},point];
				\node (1) at (2,0) [label=above:{\scriptsize $G_2$},point];
				\node (2) at (3,0) {$\cdots$};
				\node (3) at (4,0) [label=above:{\scriptsize $H_{\infty}$},point];
				\node (4) at (6,0) [label=above:{\scriptsize $H_{\omega}$},point];
				
				\draw[->] (0) edge node[above] {$K_1$} (1);
				\draw[->] (3) edge node[above] {$H_{\infty}$} (4);
			\end{tikzpicture}
		\end{center}
		
		By \cite[Section 4.3]{Wilkes}, this graph of pro-$p$ groups is injective and by \cite[Section 4.4]{Wilkes} $J=\langle G_1, H_\omega\rangle$. Since $G_1$ is finite and $H_\omega$ is 2-generated, $J$ is finitely generated (in fact for $p=2$ the group $J$ is 3-generated).  Collapsing the right edge we shall get the reduced graph of pro-$p$ groups  since no vertex group  equals to an edge group of an incident edge. Note that the latter graph of groups has a unique vertex $\infty$ whose vertex group is infinite and isomorphic to $\F_p\wr \Z_p$ which does not split over a finite $p$-group, so satifies the hypotheses of Theorem \ref{a virtual action}.

	\end{exa}

	\bigskip
	{\it Author's Adresses:}
	
	\medskip
	Mattheus P. S. Aguiar\\
	Departamento de Matem\'atica,\\
	~Universidade de Bras\'\i lia,\\
	70910-900 Bras\'\i lia DF,\\
	Brazil
	
	mattheus@mat.unb.br
	
	\medskip
	Pavel A. Zalesski\\
	Departamento de Matem\'atica,\\
	~Universidade de Bras\'\i lia,\\
	70910-900 Bras\'\i lia DF\\
	Brazil
	
	pz@mat.unb.br

\end{document}